\documentclass[12pt]{article}

\usepackage{amsthm, amsmath, amssymb, amsrefs}
\usepackage[dvipdfmx]{graphicx}

\numberwithin{equation}{section}

\theoremstyle{plain}	     
\newtheorem{thm}{Theorem}[section] 
\newtheorem{cor}[thm]{Corollary}

\theoremstyle{definition}

\theoremstyle{remark}

\newcommand{\disp}{\displaystyle}

\begin{document}

\title{Redheffer-type inequalities
for generalized trigonometric functions 
\footnote{The work of S. Takeuchi was supported by JSPS KAKENHI Grant Number 17K05336.}}
\author{Shimpei Ozawa and Shingo Takeuchi \\
Department of Mathematical Sciences\\
Shibaura Institute of Technology
\thanks{307 Fukasaku, Minuma-ku,
Saitama-shi, Saitama 337-8570, Japan. \endgraf
{\it E-mail address\/}: shingo@shibaura-it.ac.jp \endgraf
{\it 2010 Mathematics Subject Classification.} 
26D05, 33E05}}

%\date{}

\maketitle

\begin{abstract}
Famous Redheffer's inequality is generalized to a class of 
anti-periodic functions. We apply the novel inequality to
the generalized trigonometric functions and establish
several Redheffer-type inequalities for these functions.
\end{abstract}

\textbf{Keywords:} 
Redheffer's inequality,
Generalized trigonometric functions,
$p$-Laplacian

%%%%%%%%%%%%%%%%%%%%%%%%%%%%%%%%%%%%%%%%%%%%%%%%%%

\section{Introduction}

Redheffer's inequality
$$\frac{\pi^2-x^2}{\pi^2+x^2} \leq \frac{\sin{x}}{x}, \quad x \neq 0,$$
was proposed by R. Redheffer in 1968 as an advanced problem
in Amer. Math. Monthly \cite[No. 5642]{JCWSRRZ1968}.
In another issue of the journal \cite[No. 5642]{RW1969}, 
the inequality was titled \textit{A delightful 
inequality} and J.P. Williams gave a proof relying on the infinite product 
representation of the sine function.
Since then, Redheffer's inequality has been widely studied in the area of inequality
(see \cite{SB2015} and the references given there).

In 2015, S\'{a}ndor and Bhayo \cite{SB2015} presented a novel interesting proof
of Redheffer's inequality, which is based on the elementary calculus.
Their proof has the potential to generalize Redheffer's inequality 
to more functions than just the sine function,
as well to the cosine function, the hyperbolic sine and cosine functions
(\cite{CZQ2004}), and to Bessel functions (\cite{B2007}).

In this paper, inspired by the proof of \cite{SB2015}, 
we will generalize Redheffer's inequality so that it can be applied 
to a class of anti-periodic functions including the sine function.
To be precise, we establish a Redheffer-type inequality for a 
function $S$ that satisfies the following conditions:
there exist an $a \in (0,\infty)$ and a finite subset $P \subset (0,a)$ such that
\begin{enumerate}
\item[(S1)] $S(-x)=-S(x)$ and $S(a+x)=-S(x)$ for $x \in [0,\infty)$;
\item[(S2)] $0<S(x)<x$ for $x \in (0,a)$;
\item[(S3)] $S \in C([0,a]) \cap C^1([0,a)) \cap C^2([0,a) \setminus P)$;
\item[(S4)] $S'(x)^2-S''(x)S(x) \geq 1$ for $x \in [0,a) \setminus P$.
\end{enumerate}

It is clear that $S(x)$ is odd, continuous, piecewise smooth 
in $\mathbb{R}$ and 
anti-periodic with period $a$ in $[0,\infty)$; $S(na)=0$, 
$(-1)^nS(x)>0$ for $x \in (na,(n+1)a)$ and $n \in \mathbb{Z}$;
and $S(x)<x$ for $x \in (0,\infty)$ 
and $S \in C^1((-a,a))$ with $S'(0)=1$.

%\begin{lem}
%\label{lem:S}
%$S$ is an anti-periodic $C^1$-function with antiperiod $b$ in $\mathbb{R}$ satisfying $S(0)=S(b)=0$, $S(x)>0$ for $x \in (0,b)$ and $S(x)<x$ for $x>0$.
%\end{lem}

We begin with a general result on such a function $S$.

\begin{thm}
\label{thm:GRI}
Let $S$ be a function satisfying the conditions {\rm (S1)-(S4)}. Then,
\begin{equation}
\label{eq:gri}
\frac{a^2-x^2}{a^2+x^2} \leq \frac{S(x)}{x}, \quad x \neq 0.
\end{equation}
\end{thm}

It is worth pointing out that Redheffer's inequality follows immediately from 
Theorem \ref{thm:GRI}, since $S(x)=\sin{x}$ satisfies the conditions
(S1)-(S4) with $a=\pi$ and $P=\emptyset$, especially
$S'(x)^2-S''(x)S(x) \equiv 1$. 

Theorem \ref{thm:GRI} yields Redheffer-type inequalities for generalized 
functions of the trigonometric sine and cosine
functions. Before stating the inequalities, we will define 
generalized trigonometric functions.

Let $1<p,\ q<\infty$ and 
$$F_{p,q}(x):=\int_0^x \frac{dt}{(1-t^q)^{1/p}}, \quad x \in [0,1].$$
We will denote by $\sin_{p,q}$ the inverse function of $F_{p,q}$, i.e.,
$$\sin_{p,q}{x}:=F_{p,q}^{-1}(x).$$
Clealy, $\sin_{p,q}{x}$ is an increasing function in $[0,\pi_{p,q}/2]$ to $[0,1]$,
where
$$\pi_{p,q}:=2F_{p,q}(1)=2\int_0^1 \frac{dt}{(1-t^q)^{1/p}}.$$
We extend $\sin_{p,q}{x}$ to $(\pi_{p,q}/2,\pi_{p,q}]$ by $\sin_{p,q}{(\pi_{p,q}-x)}$
and to the whole real line $\mathbb{R}$ as the odd $2\pi_{p,q}$-periodic 
continuation of the function. 
Since $\sin_{p,q}{x} \in C^1(\mathbb{R})$,
we also define $\cos_{p,q}{x}$ by $\cos_{p,q}{x}:=(\sin_{p,q}{x})'$,
where ${}':=d/dx$.
Then, it follows that 
$$|\cos_{p,q}{x}|^p+|\sin_{p,q}{x}|^q=1.$$
In case $(p,q)=(2,2)$, it is obvious that $\sin_{p,q}{x},\ \cos_{p,q}{x}$ 
and $\pi_{p,q}$ are reduced to the ordinary $\sin{x},\ \cos{x}$ and $\pi$,
respectively. 
This is a reason why these functions and the constant are called
\textit{generalized trigonometric functions} (with parameter $(p,q)$)
and the \textit{generalized $\pi$}, respectively. 

The generalized trigonometric functions are well studied in the context of 
nonlinear differential equations (see \cite{KT2019} and the references given there). 
Suppose that $u$ is a solution of 
the initial value problem of $p$-Laplacian
$$-(|u'|^{p-2}u')'=\frac{(p-1)q}{p} |u|^{q-2}u, \quad u(0)=0,\ u'(0)=1,$$
which is reduced to the equation $-u''=u$ of simple harmonic motion for
$u=\sin{x}$ in case $(p,q)=(2,2)$.  
Then, 
$$(|u'|^p+|u|^q)'=\left(\frac{p}{p-1}(|u'|^{p-2}u')'+q|u|^{q-2}u\right)u'=0.$$
Therefore, $|u'|^p+|u|^q=1$, hence it is reasonable to define $u$ as 
a generalized sine function and $u'$ as a generalized cosine function.
Indeed, it is possible to show that $u$ coincides with $\sin_{p,q}$ defined as above.
The generalized trigonometric functions are often applied to
the eigenvalue problem of $p$-Laplacian.

Now, to the authors' knowledge, no Redheffer-type inequalities have 
been obtained for the generalized trigonometric functions.
Applying Theorem \ref{thm:GRI} to $\sin_{p,q}{x}$, we can prove 
the following inequalities.

\begin{thm}
\label{thm:GRI2}
Let $2 \leq p,\ q <\infty$. Then,
\begin{equation}
\label{eq:GRI}
\frac{\pi_{p,q}^2-x^2}{\pi_{p,q}^2+x^2} \leq \frac{\sin_{p,q}{x}}{x}, \quad x \neq 0.
\end{equation}
In particular, for $2 \leq p <\infty$,
\begin{equation*}
\frac{\pi_{p}^2-x^2}{\pi_{p}^2+x^2} \leq \frac{\sin_{p}{x}}{x}, \quad x \neq 0,
\end{equation*}
where $\sin_p{x}:=\sin_{p,p}{x}$ and $\pi_p:=\pi_{p,p}$.
\end{thm}

Regarding the cosine function,
Chen, Zhao and Qi \cite{CZQ2004} prove the Redheffer-type inequality
\begin{equation}
\label{eq:cos}
\frac{\pi^2-x^2}{\pi^2+x^2} <\cos{\frac{x}{2}}, \quad x \in (0,\pi).
\end{equation}
Theorem \ref{thm:GRI2}
yields the following generalization of \eqref{eq:cos} and 
also offers an alternative proof of \eqref{eq:cos} than that given in \cite{CZQ2004}.

\begin{cor}
\label{cor:cos}
Let $2 \leq q <\infty$. Then,
\begin{equation*}
\frac{\pi_{q^*,q}^2-x^2}{\pi_{q^*,q}^2+x^2}<\cos_{q^*,q}^{q^*-1}{\frac{x}{2}}, 
\quad x \in (0,\pi_{q^*,q}),
\end{equation*}
where $q^*:=q/(q-1)$.
\end{cor}

\section{Proofs of results}

In this section, we give the proofs of Theorems \ref{thm:GRI}, \ref{thm:GRI2} and
Corollary \ref{cor:cos}.

\subsection{Proof of Theorem \ref{thm:GRI}}

We follow the idea of S\'{a}ndor and Bhayo \cite{SB2015}.

For $x=\pm a$, the equality of \eqref{eq:gri} clearly holds.
Since the both sides of \eqref{eq:gri}
are even functions by (S1), it is sufficient to prove 
the strict inequality of \eqref{eq:gri} for $x \in (0,a) \cup (a,\infty)$.

It is easy to show
\begin{equation}
\label{eq:>2K}
\frac{a^2-x^2}{a^2+x^2} < \frac{S(x)}{x}, \quad x \in (a,\infty).
\end{equation}
Indeed, there exists $t>0$ such that $x=a+t$.
Using the anti-periodicity of $S$ in (S1), we know
\begin{align*}
\frac{S(x)}{x}-\frac{a^2-x^2}{a^2+x^2}
&=-\frac{S(t)}{a+t}+\frac{2at+t^2}{2a^2+2at+t^2}\\
&=\frac{t}{a+t} \left( \frac{2a^2+3at+t^2}{2a^2+2at+t^2}-\frac{S(t)}{t} \right).
\end{align*}
Since $S(t)/t<1$ for all $t \in (0,\infty)$, we have \eqref{eq:>2K}.

It remains to show that
\begin{equation}
\label{eq:0K}
\frac{a^2-x^2}{a^2+x^2} <\frac{S(x)}{x}, \quad x \in (0,a).
\end{equation}
From (S2), 
%\eqref{eq:0K} 
the inequality above is equivalent to the following inequality.
\begin{equation}
\label{eq:main}
a^2 < f(x), \quad x \in (0,a),
\end{equation}
where
\begin{equation}
\label{eq:fx}
f(x)=\frac{x^2(x+S(x))}{x-S(x)}, \quad x \in (0,a).
\end{equation}

Now, we will prove \eqref{eq:main}. Let $x \in (0,a)$.
An easy calculation yields
\begin{equation}
\label{eq:bibun}
-\frac{(x-S(x))^2}{2xS(x)}f'(x)=g(x),
\end{equation} 
where
\begin{equation}
\label{eq:g}
g(x)=x+S(x)-\frac{x^2(1+S'(x))}{S(x)}.
\end{equation}
%It follow from (ii) that
%$$S'+(S')^2-S''S \geq S'+1 >0.$$
%Therefore,
%\begin{multline*}
%S^2 g'(x)=(S'+(S')^2-S''S)\left(x-\frac{S(1+S')}{S'+(S')^2-S''S}\right)^2\\
%+S^2(1+S')\frac{-1+(S')^2-S''S}{S'+(S')^2-S''S}.
%\end{multline*}
Let $b$ be any number in $(0,a)$.
If $1+S'(b) \leq 0$, then \eqref{eq:g} yields $g(b)>0$.
We consider the case where $1+S'(b)>0$. 
First, we suppose that $1+S'(x)>0$ for $x \in [0,b]$.
Then, for $x \in (0,b] \setminus P$,
\begin{align}
S^2 g'(x)
&=((S')^2-S''S+S')x^2-2S(1+S')x+S^2(1+S') \notag \\
&=(1+S')(x-S)^2+((S')^2-S''S-1)x^2, \label{eq:g'}
\end{align}
where $S=S(x)$.
From (S4), the right-hand side of \eqref{eq:g'} is positive in $(0,b] \setminus P$.
Hence, $g'(x)>0$ in $(0,b] \setminus P$. Since $g \in C((0,b])$ and $S'(0)=1$, 
$g(x)$ is strictly increasing in $(0,b]$ and 
$$g(b)>\lim_{x \to +0}g(x)=0.$$
Next we suppose that there exists $c \in (0,b)$ such that
$1+S'(x)>0$ for $x \in (c,b]$ and $1+S'(c)=0$.
Then, in a similar way as above, we can see that
$g(x)$ is strictly increasing in $(c,b]$ and
$$g(b)>\lim_{x \to c+0}g(x)=c+S(c)>0.$$
In either case, we obtain $g(b)>0$.
Thus, using \eqref{eq:bibun}, we have $f'(x)<0$, 
hence $f(x)$ is strictly decreasing in $(0,a)$.
Since $f \in C((0,a])$,
$$f(x)>f(a)=a^2.$$
This is the desired conclusion \eqref{eq:main}.
\qed

%\begin{rem}
%Though the condition (iv) is not always necessary 
%to prove \eqref{eq:main}, it is a sufficient condition for $f(x)$ to be
%strictly decreasing and can be checked easily.
%\end{rem}

\subsection{Proof of Theorem \ref{thm:GRI2}}
Let $2 \leq p,\ q<\infty$. Then, we can see that $S(x)=\sin_{p,q}{x}$ satisfies
the conditions (S1)-(S4) with $a=\pi_{p,q}$ and $P=\{\pi_{p,q}/2\}$. 
Indeed, (S1) and (S2) are easily checked. 
For (S3),  it is known that 
$\sin_{p,q}{x} \in C^1(\mathbb{R}) \cap C^2(\mathbb{R} \setminus Z)$,
where $Z=\{ (2n+1)\pi_{p,q}/2\ |\ n \in \mathbb{Z}\}$ (\cite[Proposition 2.1]{O1984}).
Finally, (S4) is proved as follows. For $x \in [0,\pi_{p,q}/2) \cup (\pi_{p,q}/2,\pi_{p,q})$,
\begin{equation}
\label{eq:double}
S''(x)=-\frac{q}{p}\sin_{p,q}^{q-1}{x}|\cos_{p,q}{x}|^{2-p}
\end{equation}
and
\begin{align*}
S'(x)^2-S''(x)S(x)-1
&=\cos_{p,q}^2{x}+\frac{q}{p}\sin_{p,q}^q{x}|\cos_{p,q}{x}|^{2-p}-1\\
&=\left(1-\frac{q}{p}\right)\cos_{p,q}^2{x}
+\frac{q}{p}|\cos_{p,q}{x}|^{2-p}-1.
\end{align*}
Since $p,\ q \geq 2$, it is easy to show that $h(t)=(1-q/p)t^2+(q/p)t^{2-p}-1$
is nonincreasing in $(0,1]$; hence $h(t) \geq h(1)=0$ in $(0,1]$. Therefore, $S$ satisfies (S4).
Thus, we can apply Theorem \ref{thm:GRI} to $S(x)=\sin_{p,q}{x}$,
and the proof is complete. 
\qed

\subsection{Proof of Corollary \ref{cor:cos}}
Letting $p=2,\ 2 \leq q <\infty$ and $x=2^{2/q-1}y$ in \eqref{eq:GRI}, we obtain
\begin{equation}
\label{eq:lr}
\frac{\pi_{2,q}^2-(2^{2/q-1}y)^2}{\pi_{2,q}^2+(2^{2/q-1}y)^2} \leq 
\frac{\sin_{2,q}{(2^{2/q-1}y)}}{(2^{2/q-1}y)}, \quad y \neq 0.
\end{equation}
Since $\pi_{2,q}=2^{2/q-1}\pi_{q^*,q}$ by \cite[(1.10)]{T2016b}, 
the left-hand side of \eqref{eq:lr} can be rewritten as 
$$\frac{\pi_{q^*,q}^2-y^2}{\pi_{q^*,q}^2+y^2}.$$
On the other hand, we know the multiple-angle formula \cite[Theorem 1.1]{T2016b}:
for $x \in [0,\pi_{2,q}/(2^{2/q})]=[0,\pi_{q^*,q}/2]$, then
\begin{equation*}
\sin_{2,q}{(2^{2/q}x)}=2^{2/q}\sin_{q^*,q}{x}\cos_{q^*,q}^{q^*-1}{x}.
\end{equation*}
Thus, for $y \in (0,\pi_{q^*,q})$, 
the right-hand side of \eqref{eq:lr} is equal to
\begin{align*}
\frac{\sin_{q^*,q}{(y/2)} \cos_{q^*,q}^{q^*-1}{(y/2)}}{y/2}.
\end{align*}
Since $\sin_{q^*,q}{(y/2)}<y/2$ in $(0,\pi_{q^*,q})$, it is strictly less than $\cos_{q^*,q}^{q^*-1}{(y/2)}$.
This completes the proof. 
\qed

\section{Remarks}

\subsection{Estimate from above}

For \eqref{eq:0K} in the proof of Theorem \ref{thm:GRI},
it is also possible to obtain an estimate of $S(x)/x$ from above if
one supposes (S2)-(S4) and that the negative limit 
$$d:=\disp \lim_{x \to +0}\frac{S''(x)}{x}  \in (-\infty,0)$$ 
exists. Indeed, in this case,
we have seen that $f(x)$, defined as \eqref{eq:fx}, is strictly decreasing in $(0,a]$,
hence l'Hospital rule yields
$$a^2=f(a)<\lim_{x \to +0}f(x)=-\frac{12}{d}$$
and
$$\frac{a^2-x^2}{a^2+x^2} < \frac{S(x)}{x}
< \frac{12+dx^2}{12-dx^2}, \quad x \in (0,a).$$

Applying this inequality to $S(x)=\sin_{p,2}{x}$ and $a=\pi_{p,2}$, 
we see that $d=-2/p$ by \eqref{eq:double} and for $p \in [2,\infty)$,
$$\frac{\pi_{p,2}^2-x^2}{\pi_{p,2}^2+x^2} < \frac{\sin_{p,2}{x}}{x}
< \frac{6p-x^2}{6p+x^2}, \quad x \in (0,\pi_{p,2}).$$
In particular, the case of $p=2$ is due to 
S\'{a}ndor and Bhayo \cite[Theorem 1]{SB2015}.

\subsection{Powers in the L.H.S.}

\begin{figure}[htbp]
\begin{center}
\includegraphics[width=8cm,clip]{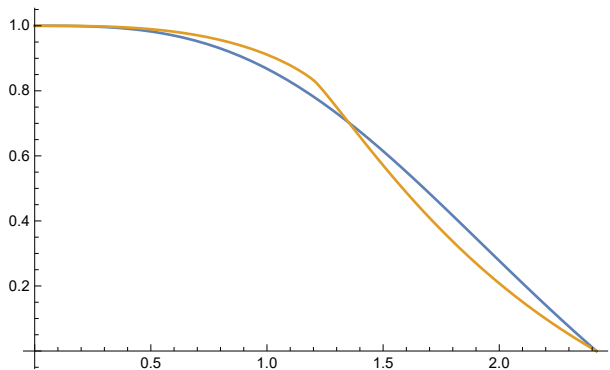}
\caption{Graphs of $(\pi_{p,q}^q-|x|^q)/(\pi_{p,q}^q+|x|^q)$ and $\sin_{p,q}{x}/x$
in $(0,\pi_{p,q}]$ for $(p,q)=(3,3)$.}
\label{fig:dame}
\end{center}
\end{figure}

We mention the left-hand side of \eqref{eq:GRI} in Theorem \ref{thm:GRI2}.
Paredes and Uchiyama \cite[Theorem 1.1]{PU2003} show that $\sin_{p,q}{x}$ has
a convergent expansion near $x=0$ as
$$\sin_{p,q}{x}=x-\frac{1}{p(q+1)}|x|^{q}x+\frac{1-p+3q-pq}{2p^2(q+1)(2q+1)}|x|^{2q}x+\cdots.$$
From the expression, we see that $\sin_{p,q}{x}/x$ can be expressed in terms of 
power series of $|x|^q$. In this sense, 
one may expect that if the left-hand side of \eqref{eq:GRI} 
is replaced with
$$\frac{\pi_{p,q}^q-|x|^q}{\pi_{p,q}^q+|x|^q},$$
it will always hold for all $p,\ q \in (1,\infty)$.
Actally, it certainly holds near $x=0$, however it does not hold near
$x=\pi_{p,q}$ (see Figure \ref{fig:dame}). 

\section*{Acknowledgement}
The authors would like to thank Professor Feng Qi for informing the paper
\cite{CZQ2004}.

\bibliographystyle{amsrefs}
\bibliography{2020redheffer}

\end{document}